\documentclass[12pt,a4paper]{amsart}
\usepackage{latexsym}
\usepackage{amssymb}

\usepackage{amssymb, verbatim}
\usepackage{graphics}  
\usepackage{graphpap}  
\usepackage{amssymb, epic, xypic} 

\makeatletter
\renewcommand{\@begintheorem}[2]{
\rm \trivlist \item [\hskip \labelsep {\bf #2\ \ #1.}]
                                }
\makeatother

\makeatletter

\DeclareFontFamily{U}{cyr}{}
\DeclareFontShape{U}{cyr}{m}{n}{
  <5> wncyr5 <6> wncyr6 <7> wncyr7 <8> wncyr8 <9> wncyr9 <10->
wncyr10}{}
\DeclareMathAlphabet{\mathcyr}{U}{cyr}{m}{n}

\input cyracc.def
\input epsf

\newcommand{\ts}{\vspace{\baselineskip}\noindent{\bf Proof.}$\;\;$}
\newcommand{\ZZ}{{\bf Z}}
\newcommand{\QQ}{{\bf Q}}

\newcommand{\CC}{{\bf C}}

\newcommand{\FF}{{\bf F}}

\newcommand{\PP}{{\bf P}}

\newcommand{\ccH}{{\mathcal H}}
\newcommand{\ccJ}{{\mathcal J}}

\newcommand{\ccS}{{\mathcal H}}

\newcommand{\bp}{\mbox{\bf p}}
\newcommand{\bu}{\mbox{\bf u}}
\newcommand{\bx}{\mbox{\bf x}}
\newcommand{\bv}{\mbox{\bf v}}
\newcommand{\by}{\mbox{\bf y}}
\newcommand{\bes}{\begin{equation*}}
\newcommand{\ees}{\end{equation*}}

\title{Some equations for the universal Kummer variety}
\author{Bert van Geemen}
\address{Dipartimento di Matematica, Universit\`a di Milano,
Via Saldini 50, 20133 Milano, Italia}

\begin{document}

\begin{abstract}
We give a method to find quartic equations for 
Kummer varieties and we give some explicit examples. 
From these equations for $g$-dimensional Kummer varieties one 
obtains equations for the moduli space of $g+1$-dimensional Kummer varieties.
These again define modular forms which vanish on the period matrices of Riemann surfaces. 
The modular forms that we find for $g=5$ appear to be
new and of lower weight than known before.
\end{abstract}

\maketitle

\section*{Introduction}
The classical Kummer surface is a quartic surface in $\PP^3$ with 16 nodes.
It is the quotient of an abelian surface by the involution $x\mapsto -x$.
The translations by points of order two on the abelian surface 
induce projective transformations that map the Kummer surface into itself.
This group action lifts to an action of a finite (non-abelian) 
Heisenberg group on $\CC^4$. 

Using classical theta function theory, the equation of this surface can be written as a Heisenberg invariant quartic polynomial whose coefficients are polynomials of degree $12$ in the second order theta constants. Thus this polynomial is best seen as a polynomial of (bi)degree $(12,4)$ in $\CC[\bu,\bx]$ which,
upon substituting the second order theta constants
for the variables $\bu$,
gives the equation of the Kummer surface.
Such a polynomial is called an equation for the universal Kummer variety.

More generally, the Kummer variety of a (principally polarized) abelian variety of dimension  $g$ admits a map
to a projective space of dimension $2^g-1$.
It is defined by second order theta functions.
For a general Kummer variety, the ideal of the image is generated by homogeneous polynomials of degree at most four.
Equations of degree four were given in \cite[Theorem 3.6]{K}.

A classical result asserts that Kummer surfaces are tangent hyperplane sections 
of the Igusa quartic threefold in $\PP^4$.
Somewhat remarkably, generalizations of the Igusa quartic
are quite useful for finding equations for 
Kummer varieties and their moduli space.
In  Proposition \ref{iq}
we show that a generalized Igusa equation (see Definition \ref{giq})
provides a Heisenberg invariant quartic equation for the universal 
Kummer variety.
In case of dimension $g=3,4$ we showed (with computer computations) that
there exist generalized Igusa equations of degree four. 
This implies that there are equations of degree $(12,4)$ 
for the universal Kummer variety in these cases 
(see Sections \ref{iq3}, \ref{iq4}).

After having read a first version of this paper, R.\ Salvati Manni suggested that 
our method would also apply to the non-Heisenberg invariant quartic equations.
This is indeed the case. 
A generalized Igusa equation of degree $d$, for $g$-dimensional abelian varieties, will provide a universal non-Heisenberg invariant quartic equation of bidegree 
$(4(d-1),4)$ for the  Kummer varieties of dimension $g+1$ (see Proposition \ref{liftrelk}). 
Thus we find explicit quartic, non-Heisenberg invariant, 
equations for the universal Kummer variety 
in genus $3,4,5$ of degree  $(12,4)$, see Section \ref{nhg3}.

The equations for the moduli space of Kummer varieties
are the universal Kummer equations of degree $(d,0)$
for some $d$. 
In section \ref{ems} we show that a generalized Igusa equation of degree $d$, 
for $g$-dimensional abelian varieties, provides an equation of degree $4d$ for the moduli space of $g+1$-dimensional abelian varieties.
Finally we recall that such an equation leads to an equation 
for the locus of Jacobians inside the moduli space of abelian varieties of
dimension $g+2$. 

This paper was motivated by the papers \cite{RSSS} and \cite{GSM}
where an explicit equation of the Coble quartic was given.
The singular locus of this quartic hypersurface in $\PP^7$ is
the Kummer variety of the Jacobian of a non-hyperelliptic curve. 
The Coble quartic is also the moduli space of rank two bundles 
with trivial determinant on the curve.
The equation for the Coble quartic is in fact an equation for the universal Kummer variety of degree $(28,4)$.
It is not (yet) clear how it is related to the 27 equations of bidegree
$(12,4)$ given in Section \ref{iq3} of this paper.

It is a pleasure to acknowledge the discussions I had with Q.\ Ren,
R.\ Salvati Manni and B.\ Sturmfels.

\section{Basics}

\subsection{Theta functions}
For $\tau\in\ccS_g$, the Siegel upper halfspace of $g\times g$ symmetric complex matrices with positive definite imaginary part,
one defines a (principally polarized) abelian variety $A_\tau$ by
$A_\tau:=\CC^g/(\ZZ^g+\ZZ^g\tau)$, were the vectors are row vectors. 
We denote by $L_\tau$ a symmetric line bundle on $A_\tau$ defining the principal polarization. 
 
The classical theta functions with characteristics $\epsilon,\epsilon'\in\{0,1\}^g$
are the holomorphic functions on 
$\ccS_g\times \CC^g$, defined by
$$
\theta[{}^\epsilon_{\epsilon'}](\tau,z)\,:=\,
\sum_{n\in\ZZ^g}\,\mbox{exp}\left[\pi i(n+\epsilon/2)\tau{}^t(n+\epsilon/2)+
2\pi i(n+\epsilon/2){}^t(z+\epsilon'/2)\right]
$$ 
(cf.\ \cite[(2.3)]{RSSS}). We will also write $\theta_m=\theta[{}^\epsilon_{\epsilon'}]$
with $m=[{}^\epsilon_{\epsilon'}]$. The function $\theta_m$ defines a global section of a translate of $L_\tau$ by a point of order two.

The second order theta functions are defined as:
$$
\Theta[\sigma](\tau,z)\,:=\,\theta[{}^\sigma_0](2\tau,2z)\qquad(\sigma\in\{0,1\}^g)~.
$$
They provide a basis of $H^0(A_\tau,L_\tau^{\otimes 2})$. 
The (totally symmetric) line bundle $L_\tau^{\otimes 2}$ is intrinsically defined by the principal polarization on $A_\tau$.

\subsection{The Kummer variety} The second order theta functions are all even functions, $\Theta[\sigma](\tau,z)=\Theta[\sigma](\tau,-z)$. The map 
$$
\Theta_\tau\,:\,A_\tau\,\longrightarrow\,\PP^{2^g-1},\qquad
z\,\longmapsto\,(\ldots:\Theta[\sigma](\tau,z):\ldots)~,
$$
which they define thus factors over the quotient variety $A_\tau/\{\pm 1\}$
which is known as the Kummer variety of the abelian variety $A_\tau$. 
The singular locus of the Kummer variety consists of $2^{2g}$ points which are the images of the fixed points of the map $x\mapsto -x$ on $A_\tau$, 
these are the 2-torsion points of $A_\tau$. 
In case the ppav $A_\tau$ is not a product of ppav's of lower dimension, 
the image of $\Theta_\tau$ is isomorphic to its Kummer variety.

\subsection{The universal Kummer variety}
The coefficients of the equations for the Kummer variety which we will give in this paper
depend only on the coordinates of the point $\Theta_\tau(0)$. 
This leads us (similar to \cite[Section 3]{RSSS}) to consider the map
$$
\Theta\,:\,
\ccS_g\,\times\,\CC^g\,\longrightarrow \,\PP^{2^g-1}\times\PP^{2^g-1}~,
\qquad (\tau,z)\,\longmapsto\,(u,x)\,=\,
\big(\Theta_\tau(0),\,\Theta_\tau(z)\big)~.
$$
The image of this map is a quasi-projective variety of dimension $g+g(g+1)/2$.
The closure of the image is called the universal Kummer variety ${\mathcal K}_g(2,4)$
and we will exhibit polynomials in $\CC[\bu,\bx]$  (actually in $\QQ[\bu,\bx]$)
which are in the ideal ${\mathcal I}_g$ of ${\mathcal K}_g(2,4)$.
In particular, for genus $2,3,4$ we show that there are such polynomials of bidegree $(12,4)$
(these are new for $g=3,4$). 

As $\Theta_\tau(0)$ lies in the image of the Kummer variety, any universal
Kummer equation $F(\bu,\bx)$ gives an equation $F(\bu,\bu)$ for the image 
of the Siegel upper half space $\ccS_g$ in $\PP^{2^g-1}$ under the map
$$
\ccS_g\,\longrightarrow\,\PP^{2^g-1},\qquad \tau\,\longmapsto\,\Theta_\tau(0)~.
$$
The new equations of bidegree $(12,4)$ for ${\mathcal K}_g(2,4)$ have the extra property that $F(\bu,\bu)$ is identically zero as polynomial in $\bu$ 
if the equation is Heisenberg invariant (as polynomial in $\bx$).
The non-Heisenberg invariant equations do give non-trivial equations
for the image and
we find equations for ${\mathcal K}_g(2,4)$ of bidegree 
$(16,0)$ for $g=3,4,5$ (which are new for $g=4,5$).

\subsection{The Heisenberg group}
Any $a\in A_\tau[2]$, the subgroup of two torsion points of $A_\tau$, 
defines a biholomorphic map $t_a:A_\tau\rightarrow A_\tau$, $x\mapsto x+a$.
The line bundles $L_\tau$ and $t_a^*L_\tau$ are isomorphic
and such an isomorphism is unique up to scalar multiple. 
However, one cannot choose the isomorphisms for the various $a$ 
in such a way as to obtain a action of $A_\tau[2]$ on  
$H^0(A_\tau,L_\tau^{\otimes 2})$. 
Instead one obtains a representation of a non-commutative Heisenberg group $H_g$, an extension of $A_\tau[2]\cong(\ZZ/2\ZZ)^{2g}$ by $\CC^\times$, on 
$H^0(A_\tau,L_\tau^{\otimes 2})$. 

The map $\Theta_\tau$ is equivariant for the action of $H_g$, where the action on
$\PP^{2^g-1}$ is generated by the following projective transformations
(which we refer to as sign changes and translations):
$$
x_\sigma\,\longmapsto\,(-1)^{\sigma{}^t\alpha}x_{\sigma}\,;\qquad
x_\sigma\,\longmapsto\,x_{\sigma+\beta}\qquad
\qquad(\sigma\in(\ZZ/2\ZZ)^g)~,
$$
where the $x_\sigma$ are the homogeneous coordinates on $\PP^{2^g-1}$
and $\alpha,\beta\in (\ZZ/2\ZZ)^g$ correspond to the elements 
$\alpha/2,\beta\tau/2\in \mbox{$\frac{1}{2}$}\Lambda_\tau/\Lambda_\tau=A_\tau[2]$ respectively.


\section {Equations for the Kummer variety}

\subsection{The general theory}
We will write a (homogeneous) polynomial of degree $d$ in the $2^g$ variables $x_\sigma$ as
$$
F\,:=\,\sum_\alpha\, c_\alpha x^\alpha \qquad 
(\in\CC[\bx]\,:=\,\CC[\ldots,x_\sigma,\ldots]),
$$
where 
$$
\alpha=(\ldots,\alpha_{\sigma},\ldots)\,\in\,\ZZ_{\geq 0}^{2^g},
\qquad 
|\alpha|\,:=\,\sum_\sigma \alpha_\sigma\,=\,d,
\qquad
x^\alpha\,:=\,\prod_\sigma x_\sigma^{\alpha_\sigma}~.
$$
The polynomial $F$ is an equation for the Kummer variety of $A_\tau$ if
the theta function 
$$
(\Theta_\tau^*F)(z)\,:=\,\sum_\alpha c_\alpha\prod_\sigma\Theta[\sigma](\tau,z)^{\alpha_\sigma}
$$
is zero as function of $z\in\CC^g$. So $\Theta_\tau^*$ substitutes $x_\sigma:=\Theta[\sigma](\tau,z)$. 

Since $F$ is homogeneous of degree $d$, 
the function $(\Theta_\tau^*F)(z))$ is a theta function of order $2d$, 
equivalently, it corresponds to a global section of $L^{\otimes 2d}_\tau$. Let $\vartheta_1,\ldots,\vartheta_N$, where $N=N_d$,
be theta functions which provide a basis of $H^0(A_\tau,L^{\otimes 2d}_\tau)$.
Then $(\Theta_\tau^*F)(z))=\sum_{i=1}^N a_i\vartheta_i(z)$ for certain complex numbers $a_i$. Moreover, $F$ is an equation for the Kummer variety if and only if $a_i=0$ for all $i$. Thus the $\CC$-vector space $I_d=I_d(\Theta_\tau(A_\tau))$ 
of equations of degree $d$ of the Kummer variety is 
$$
I_d\,:=\,\ker(\,\Theta_\tau^*\,:\,
\CC[\bx]_d\,\longrightarrow\,
H^0(A_\tau,L^{\otimes 2d}_\tau)\,)~,
$$
where 
$\CC[\bx]_d$ denotes the vector space of homogeneous
polynomials of degree $d$. (Equivalently, $I_d$ is the kernel of the natural 
map $S^dH^0(A_\tau,L^{\otimes 2}_\tau)\,\rightarrow 
H^0(A_\tau,L^{\otimes 2d}_\tau)$.)

Finding the equations of degree $d$ is thus a problem in linear algebra once the matrix of the $\CC$-linear map $\Theta_\tau^*$ w.r.t.\ bases of its domain and 
image are known. For degree $d=2$ this is matrix is the diagonal matrix with entries 
the $\theta[{}^\epsilon_{\epsilon'}](\tau,0)$. 
Unfortunately, for $d=3$ this matrix is not known explicitly. 
For $d=4$ the matrix was given in 
\cite[Proposition 4]{sjhe} and we recall part of the result in Section \ref{thd4}.

\subsection{The map $\Theta_\tau^*$ in degree two}\label{thd2}
we recall that for degree $d=2$ the matrix of $\Theta_\tau^*$, w.r.t.\  suitable bases is a diagonal matrix.
For $\epsilon,\epsilon'\in(\ZZ/2\ZZ)^g$ we define a polynomial in
$\CC[\bx]_2$ by:
$$
Q[{}^\epsilon_{\epsilon'}]\,:=\,
\sum_\sigma (-1)^{\sigma{}^t\epsilon'}\,x_\sigma x_{\sigma+\epsilon}~
$$ 
with summation over $\sigma\in(\ZZ/2\ZZ)^g$.
In case $\epsilon{}^t\epsilon'=0$ this polynomial is not identically zero. 
The $2^{g-1}(2^g+1)$ polynomials one finds in this way are a basis of 
$\CC[\bx]_2$.
(They are also eigenvectors for the action of the Heisenberg group.)
A classical theta function formula shows that (\cite[(2.8)]{RSSS})
$$
\Theta_\tau^*(Q[{}^\epsilon_{\epsilon'}])\,=\,
Q[{}^\epsilon_{\epsilon'}](\ldots,\Theta[{}{\sigma}](\tau,z),\ldots)\,=\,
\theta[{}^\epsilon_{\epsilon'}](\tau,0)\theta[{}^\epsilon_{\epsilon'}](\tau,2z)~.
$$
The functions $\theta[{}^\epsilon_{\epsilon'}](\tau,2z)$, with $\epsilon{}^t\epsilon'=0$,
are a basis of the subspace of even theta functions in 
$H^0(A_\tau,L_\tau^{\otimes 4})$. 
Therefore the matrix of $\Theta_\tau^*$ is now a diagonal matrix with 
entries the even theta constants $\theta[{}^\epsilon_{\epsilon'}](\tau,0)$.
In particular, there are quadratic equations for the Kummer variety of $A_\tau$ 
if and only if at least one of these theta constants is zero.

\subsection{The Heisenberg invariant quartics}\label{hiq}
The vector space $\CC[\bx]_4$ has the obvious basis given by the monomials $\prod_\sigma x_\sigma^{\alpha_\sigma}$ with $|\alpha|=4$.
As the map $\Theta_\tau^*$ is equivariant w.r.t.\ the action of the Heisenberg group, the matrix will be in block form if we choose bases adapted to this action.

In \cite[Proposition 1(iii)]{sjhe} it is shown that
$$
\CC[\bx]_4\,=\,\oplus_\chi\;\CC[\bx]_{4,\chi},
$$
where $\chi:H_g\rightarrow\{\pm1\}$  runs over the characters of
$H_g$, these maps factor over $(\ZZ/2\ZZ)^{2g}$.
The dimension of the eigenspace with character $\chi$ is
$(2^g+1)(2^{g-1}+1)/3$ if $\chi=0$ (i.e.\ $\chi(x)=1$ for all $x$) and it is
$(2^{g-1}+1)(2^{g-2}+1)/3$ for the $2^{2g}-1$ non-trivial characters.

A basis for $\CC[\bx]_{4,0}$, the quartic Heisenberg invariants, can be found as follows. 
Let $T\subset (\ZZ/2\ZZ)^g$ be a subgroup of order at most $4$.
So we have the following possibilities for $T$:
$$
T\,=\,\{0\},\quad\{0,\,\alpha\},\quad \{0,\alpha,\beta,\alpha+\beta\}
\qquad (\alpha,\beta\in (\ZZ/2\ZZ)^g)\qquad (\alpha\,\neq\,\beta)~.
$$
In all cases we can in fact write $T=\{0,\alpha,\beta,\alpha+\beta\}$, possibly with $\alpha=\beta=0$ or $\beta=0$. 
There are $1+(2^g-1)+(2^g-1)(2^g-2)/3$
such subgroups.
For each such subgroup $T$ the monomial $x_0x_\alpha x_\beta x_{\alpha+\beta}$
is invariant under the sign changes in the Heisenberg group. Now we simply take the sum of the monomials in the orbit of this monomial under the translations in the Heisenberg group. Then we obtain the polynomial
$$
P_T\,:=\,\sum_{\rho\in(\ZZ/2\ZZ)^g}\,
x_\rho x_{\rho+\alpha} x_{\rho+\beta} x_{\rho+\alpha+\beta}~,
\qquad T\,=\,\{0,\alpha,\beta,\alpha+\beta\}~.
$$
These $P_T$'s are a basis of $\CC[\bx]_{4,0}$.

We will need to know the partial derivative of $P_T$ w.r.t.\ $x_\sigma$. 
In case $T=\{0\}$, we have $P_T=\sum_\rho x_\rho^4$ and 
$\partial P_T/\partial x_\sigma=4x_\sigma^3$. 
In case $T=\{0,\alpha\}$, 
we have to consider the summands of $P_T$ with $\rho=\sigma$, 
which is $x_\sigma^2x_{\sigma+\alpha}^2$, 
and $\rho=\sigma+\alpha$, which is $x_{\sigma+\alpha}^2x_\sigma^2$, so 
$\partial P_T/\partial x_\sigma=4x_\sigma x_{\sigma+\alpha}^2$. 
In case $\sharp T=4$, there are four monomials which contribute to the partial derivative, and each gives $x_{\sigma+\alpha}x_{\sigma+\beta}x_{\sigma+\alpha+\beta}$.  
So we find the following result, which holds for all subgroups $T$:
$$
\frac{\partial P_T}{\partial x_\sigma}\,=\,
4x_{\sigma+\alpha}x_{\sigma+\beta}x_{\sigma+\alpha+\beta}~,
\qquad T\,=\,\{0,\alpha,\beta,\alpha+\beta\}~.
$$

\subsection{The map $\Theta_\tau^*$ in degree four}\label{thd4}
The decomposition of the image of $\Theta_\tau^*$ 
into eigenspaces for the Heisenberg group is as follows:
$$
H^0(A_\tau,L^{\otimes 8}_\tau)\,=\,
\oplus_\chi\,H^0(A_\tau,L^{\otimes 8}_\tau)_\chi~,
$$
where $\dim H^0(A_\tau,L^{\otimes 8}_\tau)_\chi$ is $2^g$
(but if $\chi\neq 0$, the subspace of even theta functions has dimension $2^{g-1}$). 
With these decompositions, we have, for all characters $\chi$ of $H_g$:
$$
\Theta_\tau^*(\CC[\ldots,x_\sigma,\ldots]_{4,\chi})\;\subset\;
H^0(A_\tau,L^{\otimes 8}_\tau)_\chi~.
$$

The `multiplication by two' map $[2]$ on the abelian variety $A_\tau$ 
(so $[2](x):=2x$ for $x\in A_\tau$) gives the following isomorphism:
$$
H^0(A_\tau,L^{\otimes 8}_\tau)_0\,=\,[2]^*H^0(A_\tau,L^{\otimes 2}_\tau),
\qquad\mbox{so the}\quad \Theta[\sigma](\tau,2z),\quad(\sigma\in(\ZZ/2\ZZ)^g)~,
$$
are a basis of $H^0(A_\tau,L^{\otimes 8}_\tau)_0$. 

We now consider the multiplication map for the case $\chi=0$:
$$
\Theta_\tau^*:\;\CC[\ldots,x_\sigma,\ldots]_{4,0}\,=\,\oplus_T\CC P_T\,
\longrightarrow\,
H^0(A_\tau,L^{\otimes 8}_\tau)_0\,=\,\oplus_\sigma \CC \Theta[\sigma](\tau,2z)~.
$$
So, given a subgroup $T=\{0,\alpha,\beta,\alpha+\beta\}$, we want to find the complex numbers $a_{\sigma,T}$ such that 
\begin{equation*}\label{mul1}
\Theta_\tau^*(P_T)\,=\,\sum_\sigma\, a_{\sigma,T}\,\Theta[\sigma](\tau,2z)
~.
\tag{\ref{thd4}.a}
\end{equation*}
Riemann's theta formula implies that (see \cite[Proposition 4]{sjhe}):
\begin{equation*}\label{mul2}
a_{\sigma,T}\,=\,\Theta[{\sigma+\alpha}] (\tau,0)
\Theta[\sigma+\beta](\tau,0) \Theta[\sigma+\alpha+\beta](\tau,0)~.
\tag{\ref{thd4}.b}
\end{equation*}
Using the derivatives of the $P_T$, there is an attractive way of to write this result:
$$
4a_{\sigma,T}\,=\,
\frac{\partial P_T}{\partial x_\sigma}(\Theta_\tau(0))~,
$$
that is, the entries of the matrix, up to a factor 4 which does not affect the kernel, are just the partial derivatives of the $P_T$'s evaluated in the point $\Theta_\tau(0)$ 
which has coordinates $x_\sigma=\Theta[\sigma](\tau,0)$. 
In other words, the matrix $(a_{\sigma,T})$ is the transposed of the 
Jacobi matrix of the polynomials $P_T$, evaluated in $\Theta_\tau(0)$.

Using Cramer's rule, applied to a $2^g\times (2^g+1)$ submatrix
of rank $2^g$, one obtains quartic Heisenberg equations for the Kummer
variety whose coefficients are of degree $3\cdot 2^g$ in the $\Theta[\sigma](\tau,0)$. 
The quartic equations given by Khaled (\cite[Theorem 3.6(b)]{K}) appear to be of degree $(2^{8g},4)$ (see the definition of the $\tilde{q}$'s on p.208 of \cite{K}), but he also considers more general embeddings of Kummer varieties.

In Section \ref{mul4g2} we recall that the classical equation for the Kummer surface in $\PP^3$ is obtained in this way. In Section \ref{mul4g3} we recall a result from \cite{RSSS} which shows that there equations of lower degree, 
$(16,4)$ rather than $(24,4)$, in case $g=3$. 
However, for $g=3,4$ there also exist equations of degree $(12,4)$ 
as we will show in Sections \ref{iq3}, \ref{iq4}.

There is a similar result for the matrix of $
\Theta_\tau^*:\CC[\ldots,x_\sigma,\ldots]_{4,\chi}\rightarrow
H^0(A_\tau,L^{\otimes 8}_\tau)_\chi$ for non-trivial $\chi$.
However, the entries of this matrix are no longer polynomials in  the $\Theta[\sigma](\tau,0)$, but they involve more general theta constants,
see Section \ref{thd4a}.

Before discussing the cases $g=2,3$, we give a characterization of the quartic, Heisenberg invariant, equations of the Kummer variety which will be exploited in the next section.

\subsection{Proposition}\label{prop1}
A quartic, Heisenberg invariant, polynomial $P$ is an equation for the Kummer variety of $A_\tau$ if and only if the point 
$\Theta_\tau(0):=(\ldots:\Theta[\sigma](\tau,0):\ldots)$ is a singular point of the 
algebraic variety defined by $P=0$ in $\PP^{2^g-1}$.

\ts Since $P\in\CC[\bx]_{4,0}$, there are $c_T\in\CC$ such that 
$P=\sum_T c_TP_T$. Using equation \ref{mul1}, the image of $P$ in $H^0(A_\tau,L_\tau^{\otimes 8})$ is
$$
\Theta_\tau^*(P)\,=\,\sum_Tc_T\Theta_\tau^*(P_T)\,=\,
\sum_{T} c_T\left(\sum_\sigma a_{\sigma,T}\Theta[\sigma](\tau,2z)\right)\,=\,
\sum_\sigma\left(\sum_Ta_{\sigma,T}c_T\right)\Theta[\sigma](\tau,2z)~.
$$
As $P$ is an equation for the Kummer variety if and only if $\Theta_\tau^*(P)=0$ on $\CC^g$, we find, using the determination of the $a_{\sigma,T}$ above, that 
this is the case exactly when
$$
0\,=\,
\sum_Ta_{\sigma,T}c_T
\,=\,
\sum_T  c_T
\frac{\partial P_T}{\partial x_\sigma}(\Theta_\tau(0))
\,=\,
\frac{\partial P}{\partial x_\sigma}(\Theta_\tau(0))
\qquad
\mbox{for all }\sigma\in(\ZZ/2\ZZ)^g~.
$$
This again is equivalent to $\Theta_\tau(0)$ being a singular point on the variety (better: subscheme) defined by $P=0$.
\qed

\subsection{The case $g=2$}\label{mul4g2}
In this case the multiplication map $\mu_{4,0}$ is given by the $4\times 5$ matrix
of derivatives of the polynomials $P_0=x_{00}^4+\ldots+x_{11}^4$, 
$$
P_1\,=\,2(x_{00}^2x_{01}^2+x_{10}^2x_{11}^2),\quad
P_2\,=\,2(x_{00}^2x_{10}^2+x_{01}^2x_{11}^2),\quad
P_3\,=\,2(x_{00}^2x_{11}^2+x_{01}^2x_{10}^2)~,
$$
and $P_{12}=4x_{00}x_{01}x_{10}x_{11}$, 
evaluated in $x_\sigma=u_\sigma$ with $u_\sigma:=\Theta[\sigma](\tau,0)$.
Using Cramer's rule, the equation of the Kummer surface 
$a_0P_0+\ldots+a_{12}P_{12}=0$ is then a determinant 
(cf.\ \cite[(1.1)]{RSSS}):
{\renewcommand{\arraystretch}{1.2}
$$
F(\bu,\bx)\,:=\,\det\left(
\begin{array}{ccccc}
P_0&P_1&P_2&P_3&P_{12}\\
u_{00}^3&u_{00}u_{01}^2&u_{00}u_{10}^2&u_{00}u_{11}^2&
u_{01}u_{10}u_{11}\\
u_{01}^3&u_{00}^2u_{01}&u_{01}u_{11}^2&u_{01}u_{10}^2&
u_{00}u_{10}u_{11}\\
u_{10}^3&u_{10}u_{11}^2&u_{00}^2u_{10}&u_{01}^2u_{10}&
u_{00}u_{01}u_{11}\\
u_{11}^3&u_{10}^2u_{11}&u_{01}^2u_{11}&u_{00}^2u_{11}&
u_{00}u_{01}u_{10}\\
\end{array}\right)\;=\,0~,
$$
} 
where we took out a factor $4$ from each of the last four rows.
The polynomial $F\in \CC[\bu,\bx]$ is of bidegree $(12,4)$,
and it is the equation defining the universal Kummer variety 
${\mathcal K}_2(2,4)$ in $\PP^3\times\PP^3$. 

\subsection{The case $g=3$}\label{mul4g3}
The basis of Heisenberg invariant quartics is given by the following $15=1+7+7$ 
polynomials:
$P_0=\sum x_\sigma^4$, $P_i=2(x_0^2x_\sigma^2+\ldots)$
where $i$ and $\sigma$ correspond through $i=4\sigma_1+2\sigma_2+\sigma_3$,
and finally $P_{ij}=4(x_0x_{\alpha}x_{\beta}x_{\alpha+\beta}+\ldots)$, where $P_{ij}=P_T$ with $T$ the subgroup of order $4$ of $(\ZZ/2\ZZ)^3$ 
generated by $\alpha,\beta$ corresponding to $i,j\in\{1,\ldots,7\}$
(and we assume that $\alpha+\beta$ corresponds to $k$ with $i<j<k$).

To be explicit, here are some of these $P_T$'s (up to $2$-power factors, they coincide with the polynomials in the $x_\sigma$ in \cite[(2.11)]{RSSS}, the others are quite similar and some can be recovered from the matrix of partial derivatives below):
{\renewcommand{\arraystretch}{1.2}
$$
\begin{array}{rcl}
P_2&=&2(x_{000}^2x_{010}^2+x_{001}^2x_{011}^2+
x_{100}^2x_{110}^2+x_{101}^2x_{111}^2),\\
P_{24}&=&4(x_{000}x_{010}x_{100}x_{110}+x_{001}x_{011}x_{101}x_{111}),\\
P_{12}&=&4(x_{000}x_{001}x_{010}x_{011}+x_{100}x_{101}x_{110}x_{111})~.
\end{array}
$$
}

In this case the multiplication matrix has size $8\times 15$, 
but it has submatrices of rank $6$. This was suggested by \cite[Lemma 8.2]{RSSS}.
We checked that any $7\times 7$ minor of the $8\times 7$ 
submatrix formed by the Jacobian matrix $(\partial P_T/\partial x_\sigma)$ of the $7$ quartics $P_2,P_3,P_{24},P_{34},P_{25},P_{35},P_{12}$ is identically zero.
This leads to an equation for the universal Kummer threefold which is a linear combination of these seven $P_T$'s.

Now we consider the following $7\times 7$ matrix, the last six rows are the partial derivatives
(divided by $4$) of the polynomials in the first row, w.r.t. the six variables $x_{001},\ldots,x_{110}$, and the derivatives are evaluated
in the point with coordinates $x_{\sigma}=u_\sigma$:

{\renewcommand{\arraystretch}{1.2}
$$
\left(
\begin{array}{ccccccc}
P_2      & P_3 &     P_{24} &     P_{34}  &    P_{25}  &    P_{35}  &    P_{12}\\
u_{001}u_{011}^2 &u_{001}u_{010}^2&u_{011}u_{101}u_{111}&
u_{010}u_{101}u_{110}& u_{011}u_{100}u_{110}&u_{010}u_{100}u_{111}& u_{000}u_{010}u_{011}\\
u_{000}^2u_{010} &u_{001}^2u_{010}&u_{000}u_{100}u_{110}& u_{001}u_{101}u_{110}&u_{000}u_{101}u_{111}&u_{001}u_{100}u_{111}& u_{000}u_{001}u_{011}\\
u_{001}^2u_{011} &u_{000}^2u_{011}&u_{001}u_{101}u_{111}&
u_{000}u_{100}u_{111}&u_{001}u_{100}u_{110}&u_{000}u_{101}u_{110}&
u_{000}u_{001}u_{010}\\
u_{100}u_{110}^2 &u_{100}u_{111}^2&u_{000}u_{010}u_{110}&
u_{000}u_{011}u_{111}&u_{001}u_{011}u_{110}&u_{001}u_{010}u_{111}&
u_{101}u_{110}u_{111}\\
u_{101}u_{111}^2 &u_{101}u_{110}^2&u_{001}u_{011}u_{111}&
u_{001}u_{010}u_{110}&u_{000}u_{010}u_{111}&u_{000}u_{011}u_{110}&
u_{100}u_{110}u_{111}\\
u_{100}^2u_{110} &u_{101}^2u_{110}&u_{000}u_{010}u_{100}&
u_{001}u_{010}u_{101}&u_{001}u_{011}u_{100}&u_{000}u_{011}u_{101}&
u_{100}u_{101}u_{111}\\
\end{array}\right)~.
$$
}

The determinant of this matrix is a linear combination $\sum_Ta_TP_T$
where each $a_T$ is a degree $18=3\cdot 6$ polynomial in the $u_\sigma$. 
Each coefficient $a_T$ turns out to be a multiple of $u_{000}u_{111}$. 
After dividing the determinant by this monomial,
we obtain a polynomial of bidegree $(16,4)$ which is an equation for the universal Kummer threefold.
Writing the polynomial as
$$
s_{010}P_2+s_{011}P_3+t_{001}P_{24}+ t_{011}P_{34}+t_{101}P_{25}
+t_{111}P_{35}+t_{100}P_{12}
$$
one obtains the coefficients $s_\sigma,t_\sigma\in\CC[\bu]$ of the polynomial in \cite[Lemma 8.2]{RSSS}.

\section{Generalized Igusa quartics}

\subsection{Singular Heisenberg invariant quartics}\label{defxsus}
In view of Proposition \ref{prop1}, we want to find quartic Heisenberg invariant 
polynomials $F$ such that all partial derivatives of $F$ vanish in 
the point $\Theta_\tau(0)$ with $\tau\in\ccS_g$. 
We consider the following variation: find quartic Heisenberg invariant
polynomials, whose coefficients are polynomials in
the coordinates of a point $u=(\ldots:u_\sigma:\ldots)\in\PP^{2^g-1}$,
such that all their partial derivatives vanish in $u$. 
So we do not require that $u=\Theta_\tau(0)$ for some $\tau$.
Such a polynomial 
$F\in\CC[\bu,\bx]$ will be an 
equation for the universal Kummer variety.

The following proposition provides a method for finding such polynomials.
It uses polynomials in variables $p_T$, 
$$
\sum_{\beta}\,d_\beta p^\beta,\qquad
p^\beta\,:=\,\prod_{T}p_{T}^{\beta_T}
$$
with $d_\beta\in\CC$ and $\beta=(\ldots,\beta_T,\ldots)$, where $T$ runs over
the set of subgroups of $(\ZZ/2\ZZ)^g$ of order at most four.
The basis $P_T=P_T(\bx)$ 
of quartic Heisenberg invariants is indexed by the same set, 
and we have an evaluation map
$$
x^*\,:\,\CC[\bp]:=\CC[\ldots,p_T,\ldots]\,\longrightarrow\,\CC[\bx],
\qquad
R\,\longmapsto\,x^*(R)\,:=\,R(\ldots,P_T(\bx),\ldots)~,
$$
which substitutes the quartic polynomial $P_T$ in the $x_\sigma$ for the variable $p_T$. Similarly, we have an evaluation map 
$$
u^*\,:\,\CC[\bp]\,\longrightarrow\,\CC[\bu],
\qquad
R\,\longmapsto\,u^*(R)\,:=\,R(\ldots,P_T(\bu),\ldots)~,
$$
so $u^*(p_T)=P_T(\bu)$.

\subsection{Definition}\label{giq}
A homogeneous polynomial $R\in\CC[\bp]$ of degree $d$ 
such that $x^*(R)=0$ is called a generalized Igusa equation of degree $d$ 
(in dimension $g$). In case $d=4$, such an $R$ will be called a generalized 
Igusa quartic (cf.\ Section \ref{iq2}).

\subsection{Definition}\label{fr}
Let $R\in \CC[\bp]_d$, 
so $R$ is a homogeneous polynomial of degree $d$ in the variables $p_T$.
Then we define 
$$
F_R\,:=\,F_R(\bu,\bx)\,=\,
\sum_T\,u^*\left(\frac{\partial R}{\partial p_T}\right)P_T(\bx)\qquad
(\,\in\,\CC[\bu,\bx]_{4(d-1),4}\,)~.
$$
For each $u\in\CC^{2^g}$, 
the polynomial $F_R(u,\bx)$ is a linear combination of the $P_T(\bx)$, 
hence it is a Heisenberg invariant quartic in $\CC[\bx]$.
The coefficient of $P_T$ is the
partial derivative of $R$ w.r.t.\ $p_T$, in which one substitutes $p_S:=P_S(\bu)$
for all indices $S$.
As ${\partial R}/{\partial p_T}$ is homogeneous of degree $d-1$, this coefficient
is a homogeneous polynomial in the $u_\sigma$ of degree $4(d-1)$. 

Upon substituting $\bx:=\bu$ in $F_R(\bx,\bu)$ one finds $R_F(\bu,\bu)=du^*R$, 
in fact:
$$
F_R(\bu,\bu)\,=\,\sum_T\,u^*\left(\frac{\partial R}{\partial p_T}\right)P_T(\bu)
\,=\, u^*\left(\sum_T\frac{\partial R}{\partial p_T}p_T\right)
\,=\,u^*(dR)~,
$$
where we used Euler's relation.

\subsection{Proposition}\label{iq}
Let $R\in \CC[\bp]_d$ be a generalized Igusa equation, so $x^*R=0$.
Then $F_R\in\CC[\bu,\bx]$ is an equation for the universal Kummer variety which moreover satisfies $F_R(\bu,\bu)=0$.

\ts
By assumption $x^*R=0$, hence also $u^*R=0$. 
The partial derivatives $\partial F_R/\partial x_\sigma$ of $F_R$, 
evaluated in $\bx=\bu$, are then 
$$
\left(\frac{\partial F_R}{\partial x_\sigma}\right)_{|\bx=\bu}\,=\,
\sum_Tu^*\left(\frac{\partial R}{\partial p_T}\right)
\left(\frac{\partial P_T}{\partial x_\sigma}\right)_{|\bx=\bu}\,=\,
\frac{\partial (u^*R)}{\partial u_\sigma}  \,=\,0~,  
$$
where we used the chain rule for differentiation.
Thus, for any $u=(\ldots,u_\sigma,\ldots)\in\CC^{2^g}$, the
Heisenberg invariant quartic $F(u,\bx)$ defines a variety in $\PP^{2^g-1}$
which is singular in the point $x=(\ldots:u_\sigma:\ldots)$.
In particular, considering the case $u=\Theta_\tau(0)$, we find that $F_R(u,\bx)$ is an equation for the Kummer variety of $A_\tau$ by Proposition \ref{prop1}.
Thus $F_R$ is an equation for the universal Kummer variety.
Finally, $F_R(\bu,\bu)=du^*R=0$.
\qed

\subsection{Remark} Proposition \ref{iq} shows that the kernel of the ring homomorphism $x^*$ provides elements of the ideal ${\mathcal I}_{g}$ of the universal Kummer variety:
$$
\ker(x^*\,:\,\CC[\bp]\,\longrightarrow\,\CC[\bx])
\,\longrightarrow\,{\mathcal I}_{g}
~,\qquad
R\,\longmapsto F_R~.
$$

The kernel of the map $x^*$ has the following geometrical interpretation.
Consider the map given by the Heisenberg invariant quartics:
$$
\PP^{2^g-1}\,\longrightarrow\,\PP^N,\qquad 
x\,=\,(\ldots:x_\sigma:\ldots)\,\longmapsto\,(\ldots:P_T(x):\ldots)~,
$$
where $N+1=(2^g+1)(2^{g-1}+1)/3$. 
For any $g$, this map factors over the Heisenberg quotient $\PP^{2^g-1}/H_g$.
In case $g=2$, one can show that it embeds $\PP^3/H_2$ into $\PP^4$ as the Igusa quartic (see also Section \ref{iq2}). (For $g\geq 3$ however, the Heisenberg invariant quartics do not generate the ring of all Heisenberg invariant
polynomials in $\CC[\bx]$). 
In any case, the homogeneous elements $R$ in the kernel of $x^*$, that is, the generalized Igusa equations, are the equations for the image of this map.  
As $2^g-1<N$ for $g\geq 2$, the kernel is non-trivial. 

A final observation is that for all $g$ this map has no base points: 
in fact, the polynomials $Q[{}^\epsilon_{\epsilon'}]$ 
from Section \ref{thd2} are a basis of the quadratic polynomials, 
hence they have no base points, 
and the $Q[{}^\epsilon_{\epsilon'}]^2$ are Heisenberg invariant quartics.

\subsection{The case $g=2$}\label{iq2}
In this case there is a unique quartic $R_2$ in the $p_T$'s which is identically zero 
as polynomial in the $x_\sigma$, i.e.\ $x^*(R_2)=0$. 
It is the well-known Igusa quartic:
$$
R_2\,:=\, p_{12}^4 + 
(p_0^2-p_1^2-p_2^2 - p_3^2)p_{12}^2 +
p_1^2 p_2^2 + p_1^2 p_3^2 + p_2^2 p_3^2 
 - 2 p_0 p_1 p_2 p_3~.
$$
From $R_2$ one finds the equation of the universal Kummer surface:
$$
F_{R_2}\,=\,u^*\Big(2p_0p_{12}^2 - 2p_1p_2p_3\Big)P_0\,+\,\ldots\,+\,
u^*\Big(4p_{12}^3+2(p_0^2-p_1^2-p_2^2 - p_3^2)p_{12}\Big)P_{12}~.
$$
In fact, $\partial R_2/\partial p_0=2(p_0p_{12}^2 - p_1p_2p_3)$ etc.
As $u^*(p_T)=P_T(\bu)$, 
the polynomial $F_{R_2}$ has bidegree $(12,4)$ in $\CC[\bu,\bx]$.
An explicit computation verifies that $F_{R_2}$ coincides with the determinant
from Section \ref{mul4g2}, up to a scalar multiple.

\subsection{The case $g=3$}\label{iq3}
The restriction of $x^*$ to polynomials of degree $\leq 3$ is injective,
but the kernel of $x^*:\CC[\bp]_4\rightarrow\CC[\bx]_{16}$ has
dimension $27$. 
An example of a quartic polynomial in the kernel is:
{\renewcommand{\arraystretch}{1.3}
$$
\begin{array}{l}
 R_3\,:=\,(p_ {14}  p_ {16}-p_1  p_{12})
 (-  p_ {24}^2 -  p_ {25}^2 +  p_ {34}^2 +  p_ {35}^2)+ p_ {14}  p_ {16} ( p_2^2 -  p_3^2)+\\
+p_ {34}  p_ {35} ( p_0  p_2 +  p_1  p_3 -  p_4  p_6 -  p_5  p_7)
- p_ {24}  p_ {25} ( p_0  p_3 +  p_1  p_2 -  p_4  p_7 -  p_5  p_6)
+\\
- p_{12} ( p_2  p_4  p_7 +  p_2  p_5  p_6 -  p_3  p_4  p_6 -  p_3  p_5  p_7)~.
\end{array}
$$
}
The corresponding basis of Heisenberg invariant quartics is given 
Section \ref{mul4g3}.
The polynomial $F_{R_3}(\bu,\bx)$ has $728$ terms.

Using the action of $Sp(6,\ZZ)$ (or the normalizer of the Heisenberg group) one can find a basis of the $27$-dimensional vector space.
As a consequence, for each $R$ in this $27$-dimensional vector space we find 
an equation of bidegree $(12,4)$ of the Kummer threefold. 
Therefore  \cite[Conjecture 8.6]{RSSS} cannot be correct as the proposed generators for the ideal of the universal Kummer threefold have bidegree $(a,b)$ with $a>12$.

\subsection{The case $g=4$}\label{iq4}
In this case we found that the kernel of $x^*$ restricted to 
$\CC[\bp]_4$ has dimension $510$.
This implies by Proposition \ref{iq} that there are equations of degree $(12,4)$ for the Kummer fourfold in $\PP^{15}$.

The space of Heisenberg invariant quartics in $2^4=16$ variables 
has dimension $51$, 
and they define a map $\PP^{15}\rightarrow\PP^{50}$.
The quartics in the kernel of $x^*$ are equations for the image of this map.
We verified the variety $Z$ (better: subscheme) in $\PP^{50}$ defined by these
equations has the image as an irreducible component.
For this we only had to find a point $x\in\PP^{15}$ such that 
tangent space $T_yZ$ has dimension $15$, where $y$ is the image of $x$.
We do not know whether $Z$ is irreducible.

\subsection{The general case}
In Section \ref{iq3}, \ref{iq4} we showed that there exist
generalized Igusa quartics for $g=3,4$. We were not able to show that such
Igusa quartics exist for any $g>2$. 
Using representation theory, 
one can make some guesses as to what to expect in general.

The action of $\Gamma_g:=Sp(2g,\ZZ)$ on the Siegel space $\ccS_g$ induces a projective projective representation of $\Gamma_g$ on the $2^g$-dimensional vector space of second order theta functions. Identifying this vector space with $\CC[\bx]_1$, the linear forms in the $2^g$ variables $x_\sigma$, 
we obtain a projective representation of $\Gamma_g$ on $\CC[\bx]$.
On the subspace of Heisenberg invariant polynomials, one obtains a 
linear representation of $Sp(2g,\FF_2)=\Gamma_g/\Gamma_g(2)$, where
$\Gamma_g(2)$ is the subgroup of matrices congruent to $I$ mod $2$.

Similarly, identifying $\CC[\bp]_1=\CC[\bx]_{4,0}$ by $p_T\mapsto P_T$,
we obtain a representation of $Sp(2g,\FF_2)$ on $\CC[\bp]$.
The representation on $\CC[\bp]_1$ is irreducible and it is the representation denoted by $\overline{V}_g$ in \cite{Frame}.
The kernel of the map $x^*:\CC[\bp]_d\,\longrightarrow\,\CC[\bx]_{4d,0}$
is a subrepresentation. This should be helpful in understanding the kernel.

The results on the kernel of $x^*$ found above fit into the following pattern.
The group $Sp(2g,\FF_2)$ has an irreducible representation of dimension 
$(2^g-1)(2^{g-1}-1)/3$ 
(on the Heisenberg invariants in $\wedge^4(\CC^{2^g})$),
denoted by $\overline{U}_g$  in \cite{Frame}, for $g\geq 2$.
The representation $Sym^2(\overline{U}_g)$ is reducible for $g>2$ 
(in case $g=2$ it is the trivial representation). 
In case $g=3$ there are two irreducible components, 
one is the trivial representation, the other has dimension $27$ 
(and this is the representation on the $27$-dimensional kernel of $x^*$).
For $g>3$, $Sym^2(\overline{U}_g)$ has a subrepresentation of dimension $2^{g-1}(2^g-1)$
(this subrepresentation is the image of the map 
$Sym^2(\overline{U}_g)\rightarrow (\wedge^8\CC^{2^g})$).
A computation shows that the complementary 
representation in $Sym^2(\overline{U}_g)$
has dimension
$$
\frac{2}{9}(2^g+1)(2^g-1)(2^{g-1}+1)(2^{g-3}-1)~
$$
for $g\geq 4$.
The dimension is thus $510$, $11594$, $210210$ for $g=4,5,6$.
We checked, using Magma, that this representation is irreducible for $g=4,5,6$.
We also found that this representation occurs, with multiplicity one,
in the representation of $Sp(2g,\FF_2)$ on the subspace
$Sym^4(\CC[\bx]_{4,0})=\CC[\bp]_4$ for $g=3,\ldots,6$. 

This suggest that $\dim\ker(x^*:\CC[\bp]_4\longrightarrow\CC[\bx]_{16})$
might have this dimension also for $g>4$.

In any case, if this kernel is non-trivial then
we would have equations of bidegree $(12,4)$ for 
the universal Kummer variety for all $g\geq 2$.
One still has to check that not all coefficients, 
that are degree $12$ polynomials in
the $u_\sigma$, are zero on the image of $\ccS_g$ in $\PP^{2^g-1}$.

\section{Non-Heisenberg invariant quartic equations}\label{nhiqe}

\subsection{The map $\Theta_\tau^*$ in degree four (bis)}\label{thd4a}
In section \ref{thd4} we considered the map $\Theta_\tau^*:\CC[\bx]_{4,0}\rightarrow H^0(A_\tau,L_\tau^{\otimes 8})_0$,
whose kernel are the quartic Heisenberg invariant equations for 
the Kummer variety of $A_\tau$.
We now consider the case of a non-trivial character $\chi$ of the Heisenberg group
$H_g$.
The action of $Sp(2g,\ZZ)$ on the non-trivial characters of $H_g$ is transitive.
Thus we concentrate on a fixed character.

It will be convenient to consider $g+1$-dimensional abelian varieties and to 
consider the character
$\chi:H_{g+1}\rightarrow\{\pm 1\}$ which is
trivial on the sign changes and which is $(-1)^{\beta_{g+1}}$ 
on the translation $v_\sigma\mapsto v_{\sigma+\beta}$
with $\beta\in (\ZZ/2\ZZ)^{g+1}$.
We introduce a polynomial ring $\CC[\bv]$ in $2^{g+1}$ variables $v_\sigma$, $\sigma\in (\ZZ/2\ZZ)^{g+1}$, which is the analogue of the ring $\CC[\bu]$
in the $2^g$ variables $u_\sigma$ for $\sigma\in (\ZZ/2\ZZ)^{g}$.
Now we define a ring homomorphism 
$$
\CC[\bp]\,\longrightarrow \,\CC[\bv],\qquad 
R\,\longmapsto\,\tilde{R}_\chi~,
$$
which is defined by the assignment
$$
p_T\,\longmapsto\,\tilde{p}_{T,\chi}\,:=\,P_T(\ldots,v_{\sigma0},\ldots)\,-\,
P_T(\ldots,v_{\sigma1},\ldots)~,
$$
where the $P_T$ are a basis of the Heisenberg invariant quartics (in the $2^g$ variables $x_\sigma$). These are thus evaluated first in $x_\sigma:=v_{\sigma0}$
and next in $x_\sigma:=v_{\sigma1}$. 

As the $P_T$ are $H_g$-invariants, 
it is easy to see that $\tilde{p}_{T,\chi}\in\CC[\bv]_{4,\chi}$.
Moreover, this homomorphism induces an isomorphism $\CC[\bp]_1\cong \CC[\bv]_{4,\chi}$.
Thus the polynomials $\tilde{p}_{T,\chi}$, where $T$ runs over the 
$(2^g+1)(2^{g-1}+1)/3$
subgroups of order at most four of $(\ZZ/2\ZZ)^g$,
form a basis of  $\CC[\bp]_1\cong \CC[\bv]_{4,\chi}$.

A basis of $H^0(A_\tau,L_\tau^{\otimes 8})_{\chi,+}$,
the even theta functions in 
$H^0(A_\tau,L_\tau^{\otimes 8})_\chi$,
is given by the $\theta[{}^{\sigma0}_{0\,1}](2\tau,4z)$, 
where $\sigma$ runs over $(\ZZ/2\ZZ)^g$ (here $\tau\in\ccS_{g+1}$, 
$z\in\CC^{g+1}$).

W.r.t.\ these bases, the map
$$
\Theta_\tau^*:\;\CC[\ldots,v_\sigma,\ldots]_{4,\chi}\,=\,
\oplus_T\CC\tilde{p}_{T,\chi}\,\longrightarrow\,
H^0(A_\tau,L_\tau^{\otimes 8})_{\chi,+}
\,=\,\oplus_\sigma \CC \theta[{}^{\sigma0}_{0\,1}](2\tau,4z)
$$
is determined by the complex numbers $b_{\sigma,\tau}$ such that
\begin{equation*}\label{mul1a}
\Theta_\tau^*(\tilde{p}_{T,\chi})\,=\,\sum_\sigma b_{\sigma,T}
\theta[{}^{\sigma0}_{0\,1}](2\tau,4z)~.
\tag{\ref{thd4a}.a}
\end{equation*}
Using Riemann's theta formula one finds (see \cite[Proposition 4]{sjhe}):
\begin{equation*}\label{mul2a}
b_{\sigma,T}\,=\,
\theta[{}^{\sigma+\alpha\,0}_{\;\;\, 0\;\;\;1}](2\tau,0)
\theta[{}^{\sigma+\beta\,0}_{\;\;\;0\;\;\; 1}](2\tau,0)
\theta[{}^{\sigma+\alpha+\beta\,0}_{\quad0\quad\; 1}](2\tau,0)~.
\tag{\ref{thd4a}.b}
\end{equation*}

\

\subsection{Proposition}\label{liftrelk}
Let $R\in\CC[\bp]_d$ be a generalized Igusa equation of degree $d$
in dimension $g$, so 
$R\in\ker(x^*:\CC[\bp]_d\rightarrow\CC[\bx]_{4d})$. 
Then the polynomial
$$
\tilde{F}_{R,\chi}\,:=\,\sum_T
\widetilde{\left(\frac{\partial R}{\partial p_T}\right)}_\chi \tilde{p}_{T,\chi}(\by)
\in\CC[\bv,\by]_{4(d-1),4}
$$
is a quartic equation for for the universal Kummer variety ${\mathcal K}_{g+1}$
which,  as polynomial in $\by$, lies in $(\CC[\bv])[\by]_{4,\chi}$.

\ts
We use the Heisenberg group $H_{g+1}$ acting on $\PP^{2^{g+1}-1}$. The coordinates on this projective space, in terms of which the action of $H_{g+1}$ is by translations and signs changes,  
will be denoted by $y_\sigma$, $\sigma\in(\ZZ/2\ZZ)^{g+1}$.
If $a\in H_{g+1}$ acts non-trivially on 
$\PP^{2^{g+1}-1}$, then it is has two eigenspaces, each isomorphic
to $\PP^{2^g-1}$. We choose one of them and call it $\PP^{2^g-1}_a$.
The subgroup of elements in $H_{g+1}$ which map this eigenspace into itself induce the action of $H_g$ on $\PP^{2^g-1}_a$.
This allows one to choose coordinates $x_\sigma$, 
$\sigma\in(\ZZ/2\ZZ)^g$,
adapted to $H_g$ on this eigenspace. 

Let $a\in H_{g+1}$ be the element 
which acts as the sign change
$$
y_{\sigma_1,\ldots,\sigma_{g+1}}\,\mapsto \,
(-1)^{\sigma_{g+1}}y_{\sigma_1\ldots\sigma_{g+1}},\quad
\mbox{then}\quad
\PP^{2^g-1}_a\,:=\,\{y\in \PP^{2^{g+1}-1}:\;
y_{\sigma_1\ldots\sigma_{g+1}}\,=\,0\;\mbox{if}\;\sigma_{g+1}=1\,\}
$$
is one of the eigenspaces of $a$ and the coordinates on $\PP^{2^g-1}_a$
are $x_{\sigma_1\ldots\sigma_g}=y_{\sigma_1\ldots\sigma_{g}0}$.

The intersection of a Kummer variety $\Theta(A_\tau)\subset\PP^{2^{g+1}-1}$ 
with the eigenspace $\PP^{2^g-1}_a$ contains $2^{2(g-1)}$ points. 
These points are the images of certain points of order four in $A_\tau$. 
In particular, the point $b:=(1/4)e_{g+1}\in\CC^{g+1}$, 
where the $e_j$ are the standard basis of $\CC^{g+1}$, 
has image $\Theta_\tau(b)\in \PP^{2^g-1}_a$. In fact one easily computes
that
$
\Theta_\tau(b)\,=\,(\ldots:y_\sigma:\ldots)$ 
with
$$
y_{\sigma_1\ldots\sigma_g1}\,=\,\Theta[\sigma_1\ldots\sigma_g1](\tau,b)
\,=\,
\theta[^{\sigma_1\ldots\sigma_g1}_{0\;\ldots\;0\;\;0}](2\tau,(1/2)e_{g+1})
\,=\,
\theta[^{\sigma_1\ldots\sigma_g1}_{0\;\ldots\;0\;\;1}](2\tau,0)
\,=\,0~,
$$
since the last theta function is odd. The other coordinates of $\Theta_\tau(b)$ are:
$$
y_{\sigma_1\ldots\sigma_g0}\,=\,\Theta[\sigma 0](\tau, b)\,=\,
\theta[{}^{\sigma 0}_{0\,0}](2\tau,(1/2)e_g)\,=\,
\theta[{}^{\sigma 0}_{0\,1}](2\tau,0)~,
$$
and we notice that these are the theta constants which appear in the 
$b_{\sigma,T}$ from equation \ref{thd4a}.b.

Given $R\in\ker(x^*:\CC[\bp]_d\rightarrow\CC[\bx]_{4d}$, 
we write the polynomial $F_R$ (cf.\ Definition \ref{fr}) as
$$
F_R\,=\,\sum_T  c_T(\bu)P_T(\bx)  \qquad\mbox{so}\quad
 c_T(\bu)\,:=\,u^*\left(\frac{\partial R}{\partial p_T}\right)~.
$$
From the proof of Proposition \ref{iq} and Section \ref{hiq} we have
$$
0\,=\,\left(\frac{\partial F_R}{\partial x_\sigma}\right)_{|\bx=\bu}
\;=\,
4\sum_T 
u_{\sigma+\alpha}u_{\sigma+\beta}u_{\sigma+\alpha+\beta}
c_T(\bu) \qquad(\in\CC[\bu])~,
$$
for all $\sigma\in(\ZZ/2\ZZ)^g$. The monomial in front of $c_T(\bu)$
is similar to $b_{\sigma,T}$ and we exploit this fact.

Let $\tau\in\ccS_{g+1}$. Then we substitute 
$$
u_\sigma\, := \,y_{\sigma_1\ldots\sigma_g0}\,=\,\Theta[\sigma 0](\tau, b) \quad  \mbox{with} 
\quad  \tau\in \ccS_{g+1},\quad  
b\,=\,e_{g+1}/4~.
$$
Then $u_{\sigma+\alpha}u_{\sigma+\beta}u_{\sigma+\alpha+\beta}$ becomes
$\theta[{}^{\sigma+\alpha\,0}_{\;\;\, 0\;\;\;1}](2\tau,0)
\theta[{}^{\sigma+\beta\,0}_{\;\;\;0\;\;\; 1}](2\tau,0)
\theta[{}^{\sigma+\alpha+\beta\,0}_{\quad0\quad\; 1}](2\tau,0)=b_{\sigma,T}$ 
and thus
$$
0\,=\,
\sum_T 
b_{\sigma,T}
c_T(\ldots,\Theta[\sigma 0](\tau, b),\ldots)\,=\,0 
$$
for all $\sigma\in(\ZZ/2\ZZ)^g$.
This again implies that the polynomial
$$
\sum_T c_T(\ldots,\Theta[\sigma 0](\tau, b),\ldots)\tilde{p}_{T,\chi}(\by)
\qquad(\in\CC[\by]_{4,\chi})
$$
is an equation for the Kummer variety of the 
$g+1$-dimensional abelian variety $A_\tau$. 

To see that this is a universal Kummer equation, 
we need to show that the coefficients are actually polynomials in the 
$\Theta[\sigma](\tau,0)$.

For a subgroup $T\subset(\ZZ/2\ZZ)^{g}$ with at most four elements,
we consider its image $\tilde{T}\subset (\ZZ/2\ZZ)^{g+1}$ under the map
$(\ZZ/2\ZZ)^{g}\hookrightarrow (\ZZ/2\ZZ)^{g+1}$, 
$(\sigma_1,\ldots,\sigma_g)\mapsto (\sigma_1,\ldots,\sigma_g,0)$.
The corresponding Heisenberg invariant quartic polynomial,
$P_{\tilde{T}}\in \CC[\bv]_{4,0}$  can be written as
$$
P_{\tilde{T}}\,=\,\sum_\sigma
v_\sigma v_{\sigma+\alpha} v_{\sigma+\beta} v_{\sigma+\alpha+\beta}
\,=\,P_T(\ldots,v_{\sigma0},\ldots)\,+\,P_T(\ldots,v_{\sigma1},\ldots)
$$
where $\sigma$ runs over $(\ZZ/2\ZZ)^{g+1}$,
$\tilde{T}=\{0,\alpha,\beta,\alpha+\beta\}$, and we split the summation over the $\sigma$ with $\sigma_{g+1}=0$ and $\sigma_{g+1}=1$
respectively.

Now we prove the following key result:
\begin{equation*}\label{key}
P_T(\ldots,\Theta[\sigma_1,\ldots\sigma_g0](\tau,b),\ldots)\,=\,
\tilde{p}_{T,\chi}(\Theta_\tau(0))~,
\tag{*}
\end{equation*}
by computing the value of $P_{\tilde{T}}(\Theta_\tau(b))$ in two different ways.

First of all, we substitute the coordinates $y_\sigma$ of 
$\Theta_\tau(b)$ in $P_{\tilde{T}}$.
As the coordinates with $\sigma_{g+1}=1$ are all zero, the result is simply
$$
P_{\tilde{T}}(\Theta_\tau(b))\,=\,P_T(\ldots,\Theta[\sigma_1,\ldots\sigma_g0](\tau,b),\ldots)~.
$$
Next we use the formulae (\ref{thd4}a,b) for $\Theta_\tau^*(P_{\tilde{T}})$:
$$
P_{\tilde{T}}(\Theta_\tau(b))\,=\,\Theta_\tau^*(P_{\tilde{T}})_{|z=b}\,=\,
\sum_\sigma 
\Theta[{\sigma+\alpha}] (\tau,0)
\Theta[\sigma+\beta](\tau,0) \Theta[\sigma+\alpha+\beta](\tau,0)
\Theta[\sigma](\tau,2b)~.
$$
The $\Theta[\sigma](\tau,2b)$ are the coordinates of the point $\Theta_\tau(2b)$.
As $a=2b$ is a 2-torsion point, 
the point $\Theta_\tau(2b)$ is obtained from the point $\Theta_\tau(0)$ by the action of the sign change of $a\in H_{g+1}$ above. In fact, one computes that 
$$
\Theta[\sigma](\tau,2b)\,=\,
\theta[{}^\sigma_0](2\tau,4b) \,=\,
\theta[{}^\sigma_0](2\tau,e_{g+1})\,=\,
(-1)^{\sigma_{g+1}}\theta[{}^\sigma_0](2\tau,0)\,=\,(-1)^{\sigma_{g+1}}
\Theta[\sigma](\tau,0)~.
$$
As $\alpha_{g+1}=\beta_{g+1}=0$
we find:
{\renewcommand{\arraystretch}{1.4}
$$
\begin{array}{rcl}
P_{\tilde{T}}(\Theta_\tau(b))&=&
\sum_{\sigma,\sigma_{g+1}=0} 
\Theta[{\sigma+\alpha}] (\tau,0)
\Theta[\sigma+\beta](\tau,0) \Theta[\sigma+\alpha+\beta](\tau,0)
\Theta[\sigma](\tau,0)\,+\,\\
&&\quad -\; 
\sum_{\sigma,\sigma_{g+1}=1} 
\Theta[{\sigma+\alpha}] (\tau,0)
\Theta[\sigma+\beta](\tau,0) \Theta[\sigma+\alpha+\beta](\tau,0)
\Theta[\sigma](\tau,0)\\
&=&
P_T(\ldots,\Theta[\sigma_1\ldots\sigma_g0](2\tau,0),\ldots)\,-\,
P_T(\ldots,\Theta[\sigma_1\ldots\sigma_g1](2\tau,0),\ldots)\\
&=&\tilde{p}_{T,\chi}(\Theta_\tau(0))~.
\end{array}
$$
}
This concludes the proof of the key result \eqref{key}.

Finally we show that the equation we found is universal.
Each partial derivative $\partial F/\partial p_T$ 
of $F$ is a polynomial in the $p_T$, 
hence the key result \eqref{key} implies that the coefficients 
are polynomials in the $\tilde{p}_{T,\chi}(\Theta_\tau(0))$, to be precise, if
$$
\frac{\partial R}{\partial p_T}\,=\,
\sum_\beta d_\beta\prod_S p_S^{\beta_S}\qquad(\in \, \CC[\bp]_{d-1})~,
$$
where $\beta=(\ldots, \beta_S,\ldots)$, and $S$ runs over the subgroups of
$(\ZZ/2\ZZ)^g$ order at most four, then 
{\renewcommand{\arraystretch}{1.5}
$$
\begin{array}{rcl}
c_T(\ldots,\Theta[\sigma 0](\tau, b),\ldots)&=&
\sum_\beta d_\beta\prod P_S(\ldots,\Theta[\sigma 0](\tau, b),\ldots)\\
&=&\sum_\beta d_\beta\prod \tilde{p}_{S,\chi}(\Theta_\tau(0))\\
&=&\widetilde{\left({\partial R}/{\partial p_T}\right)}_\chi(\Theta_\tau(0))~.
\end{array}
$$
}
Thus the coefficients of the equation of the Kummer variety are the polynomials $\widetilde{\left({\partial R}/{\partial p_T}\right)}_a\in\CC[\bv]$, which are evaluated in $\bv=\Theta_\tau(0)$. Hence we found a universal Kummer relation, which is exactly $\tilde{F}_{R,\chi}$.
\qed

\subsection{Example}\label{nhg3}
In case $g=2$, the Igusa quartic $R_2$
generates the kernel of $x^*$. 
Using the formula for $F_{R_2}$ from Section \ref{iq2}, 
one easily finds the corresponding universal Kummer equation for $g=3$: 
$\;\tilde{F}_{R_2,\chi}\,=\,$
$$
\Big(2\tilde{p}_{0,\chi}\tilde{p}_{12,\chi}^2 - 2\tilde{p}_{1,\chi}\tilde{p}_{2,\chi}\tilde{p}_{3,\chi}\Big)\tilde{p}_{0,\chi}(\by)\,+\,\ldots\,+\,
\Big(4\tilde{p}_{12,\chi}^3 + 2(\tilde{p}_{0,\chi}^2-\tilde{p}_{1,\chi}^2 -\tilde{p}_{2,\chi}^2 - \tilde{p}_{3,\chi}^2)\tilde{p}_{12,\chi}
\Big)p_{12,\chi}(\by)~.
$$
Here the basis of $\CC[\bv]_{4,\chi}$ is (cf.\ Section \ref{mul4g2}):
{\renewcommand{\arraystretch}{1.3}
$$
\begin{array}{lcl}
\tilde{p}_{0,\chi}&=&v_{000}^4+v_{010}^4+v_{100}^4+v_{110}^4-
v_{001}^4-v_{011}^4-v_{101}^4-v_{111}^4,\\
\tilde{p}_{1,\chi}&=&2(v_{000}^2v_{010}^2+v_{100}^2v_{110}^2-
v_{001}^2v_{011}^2-v_{101}^2v_{111}^2),\\
\tilde{p}_{2,\chi}&=&2(v_{000}^2v_{100}^2+v_{010}^2v_{110}^2-
v_{001}^2v_{101}^2-v_{011}^2v_{111}^2),\\
\tilde{p}_{3,\chi}&=&2(v_{000}^2v_{110}^2+v_{010}^2v_{100}^2-
v_{001}^2v_{111}^2-v_{011}^2v_{101}^2),\\
\tilde{p}_{12,\chi}&=&4(v_{000}v_{010}v_{100}v_{110}-
v_{001}v_{011}v_{101}v_{111})~.
\end{array}
$$
}

As far as I know, this equation for the universal Kummer threefold was not known before.  For each non-trivial character of $H_3$ one finds such a quartic relation 
in the corresponding eigenspace, thus we get $63$ such equations.

Similarly, using the results from Sections \ref{iq3}, \ref{iq4}, one finds non-Heisenberg invariant quartics for $g=4,5$.

\subsection{Remark} 
The quartic polynomials $Q[{}^\epsilon_{\epsilon'}]^2$ are Heisenberg invariant.
It is not hard to find their image in $\CC[\bv]$:
$$
\widetilde{(Q[{}^\epsilon_{\epsilon'}]^2)}_\chi\,=\,
Q[{}^{\epsilon\;0}_{\epsilon'0}](\bv)Q[{}^{\epsilon\;0}_{\epsilon'1}](\bv)~.
$$

\section{Equations for the moduli space}\label{ems}

\subsection{Equations for the moduli space}\label{eqms}
In this section we show that the 
generalized Igusa equations defined in Section \ref{giq} 
also produce equations for the moduli space $\Theta(\ccS_{g+1})\subset\PP^{2^{g+1}-1}$ 
(note that $g$ increases by $1$ in this process, like for the non-Heisenberg invariant equations for the universal Kummer variety in Section \ref{nhiqe}).

\subsection{Proposition}\label{liftrel}
Let $R\in\CC[\bp]_d$ be a generalized Igusa equation of degree $d$, so 
$R\in\ker(x^*:\CC[\bp]_d\rightarrow\CC[\bx]_{4d})$. 
Then the polynomial
$\tilde{R}_\chi\in\CC[\bv]_d$ is an equation for (the 
closure of) $\Theta(\ccS_{g+1})$ in  $\PP^{2^{g+1}-1}$.

\ts
Let $R\in\ker(x^*:\CC[\bp]_d\rightarrow\CC[\bx]_{4d}$).
By Proposition \ref{liftrelk} it defines the  the equation 
$\tilde{F}_{R,\chi}\in\CC[\bv,\by]$ for 
the universal Kummer variety. 
Thus $\tilde{F}_{R,\chi}(\Theta_\tau(0),\Theta_\tau(z))=0$ for all 
$\tau\in\ccS_g$ and all $z\in\CC^{g+1}$. Putting $z=0$, we see that
$\tilde{F}_{R,\chi}(\bv,\bv)$ is an equation for the moduli space. 
Euler's relation shows that
$$
\tilde{F}_{R,\chi}(\bv,\bv)\,=\, \sum_T
\widetilde{\left(\frac{\partial R}{\partial p_T}\right)}_a \tilde{p}_{T,\chi}(\bv)
\,=\,
\left(\sum_T \left(\frac{\partial R}{\partial p_T}\right)p_T\right)_{|p_T=
\tilde{p}_{T,\chi}(\bv)}
\,=\,d\widetilde{R}_\chi(\bv)~,
$$
hence also $\widetilde{R}_\chi(\bv)$ is an equation for the moduli space.
\qed

\subsection{Example $g=2$}\label{igf16}
In case $g=2$, the kernel of $x^*$ is generated by the Igusa quartic $R_2$, 
see Section \ref{iq2}. 
It is also well-known that the image of the map $\Theta:\ccS_3\rightarrow \PP^7$,
known as the Satake hypersurface,
is defined by an irreducible homogeneous polynomial $f_{16}$ of degree $16$ 
(cf.\ \cite[Proposition 3.1]{RSSS}).

The space $\CC[\bx]_{4,0}$ is five dimensional, its basis $P_T$ was given 
explicitly in Section \ref{iq2}. The $\tilde{p}_{T,\chi}$ are given in Section
\ref{nhg3}.
The polynomial $\tilde{R}_{2,\chi}$ in the eight variables $v_\sigma$ can then be computed  by substituting $p_T:=\tilde{p}_{T,\chi}$ in the polynomial $R_2$ from 
Section \ref{iq2}.
One finds that  $\tilde{R}_{2,\chi}$ is, up to multiplication by a non-zero constant, the polynomial $f_{16}$.
Thus Proposition \ref{liftrel} produces non-trivial equations.

\subsection{Example $g=3$}\label{g34mod}
We checked that the polynomial $\tilde{R}_{3,\chi}$ in the $16$ variables
obtained from the generalized Igusa equation $R_3$ from Section \ref{iq3} 
is not identically zero. It is thus an equation of degree $16$ for the image of $\ccS_4$ in $\PP^{15}$.

Classically, equations of degree $32$ 
for the moduli space were known. 
These are obtained by `rationalizing' relations between the theta constants $\theta[{}^\epsilon_{\epsilon'}](\tau,0)$ in order to obtain relations in which only
the squares $\theta[{}^\epsilon_{\epsilon'}](\tau,0)^2$ appear, 
and then express these squares in terms of the $\Theta[\sigma](\tau,0)$ 
as in Section \ref{thd2} with $z=0$. 

In \cite{FO}, Freitag and Oura show that there exists an equation of degree 24 for $\Theta(\ccS_4)$ in $\PP^{15}$.

\section{Schottky-Jung relations}

\subsection{The classical Schottky-Jung relations}\label{clasSJ}
The Schottky-Jung relations relate the theta constants 
$\theta[{}^\eta_{\eta'}](\tau,0)$, where $\tau\in\ccS_g$ is a
period matrix of the Jacobian of a genus $g$ curve $C$, to the theta constants $\theta[{}^\epsilon_{\epsilon'}](\pi,0)$, where $\pi\in\ccS_{g-1}$
is a period matrix of the Prym variety of an unramified double cover $\tilde{C}\rightarrow C$. The classical form of these relations is as follows:
\begin{equation*}\label{SJ}
\theta[{}^\epsilon_{\epsilon'}]^2(\pi,0)\,=\,
c\,\theta[{}^{\epsilon\;0}_{\epsilon'0}](\tau,0)
\theta[{}^{\epsilon\;0}_{\epsilon'1}](\tau,0)~,
\tag{SJ}
\end{equation*}
for all even characteristics $[{}^\epsilon_{\epsilon'}]$, 
where $c$ is a non-zero complex number.

Given a homogeneous polynomial in the $\theta[{}^\epsilon_{\epsilon'}]^2(\pi,0)$ 
which is identically zero as a function of $\pi\in\ccS_{g-1}$, one obtains 
from the Schottky-Jung relations a polynomial in certain
$\theta[{}^\eta_{\eta'}](\tau,0)$ and thus a holomorphic function  (actually
a modular form) on $\ccS_g$. This function is thus zero in $\tau\in\ccS_g$ 
if $\tau$ is the period matrix of a Riemann surface. 

A nice overview of the approach to the Schottky problem which uses Schottky-Jung relations and modular forms is given in \cite[Section 3]{SGr}, a quick derivation of the Schottky-Jung relations can be found in \cite[Section 6.6]{herb}.

\subsection{Example $g=4$}\label{sj4}
Let $f_{16}\in\CC[\bv]$ be the degree $16$ polynomial 
from Section \ref{igf16} which 
defines the image of $\Theta(\ccS_3)\subset\PP^7$. 
As the quadratic polynomials $Q[{}^\epsilon_{\epsilon'}]\in \CC[\bv]$
are a basis of  $ \CC[\bv]_2$, one can write $f_{16}$ as  
a degree $8$ polynomial $\bar{f}_8$
in the  $Q[{}^\epsilon_{\epsilon'}]$:
$$
f_{16}(v_{000},\ldots,v_{111})\,=\,
\bar{f}_8(\ldots,Q[{}^\epsilon_{\epsilon'}](\bv),\ldots)~.
$$
For $\pi\in\ccS_3$ one has 
$Q[{}^\epsilon_{\epsilon'}](\Theta_\pi(0))=
\theta[{}^\epsilon_{\epsilon'}]^2(\pi,0)$ 
(see Section \ref{thd2}). Thus if we define
$$
F(\tau)\,:=\,\bar{f}_8(\ldots,\theta[{}^{\epsilon\;0}_{\epsilon'0}](\tau,0)
\theta[{}^{\epsilon\;0}_{\epsilon'1}](\tau,0),\ldots)
$$
then for a period matrix of a Riemann surface $\tau\in\ccS_4$ we have
$$
F(\tau)\,=\,c^{-8}\,f_{16}(\Theta_\pi(0))\,=\,0~.
$$
Schottky verified that $F$ is not identically zero on $\ccS_4$ and
Igusa \cite{sch} (and independently Freitag \cite{Fsch}) showed that $F(\tau)=0$ implies that $\tau$ is in the Jacobi locus $\ccJ_4\subset\ccS_4$, 
which is the closure of the locus of period matrices of Riemann surfaces in $\ccS_4$.

\subsection{Example $g=5$}
Proceeding as in Example \ref{sj4},
the degree $16$ equations for $\Theta(\ccS_4)\subset\PP^{15}$
which we found in Section \ref{g34mod}
lead to polynomials of degree $16$ in the $\theta[{}^\eta_{\eta'}](\tau,0)$'s 
with $\tau\in\ccS_5$.
These are thus modular forms weight $8$ on $\Gamma_5(4,8)$ which 
are zero on the Jacobi locus $\ccJ_5$. 

From \cite[Theorem 1.3]{nick},
which generalizes a result of Grushevsky and Salvati Manni, 
we obtain that these modular forms are actually cusp forms.
In fact, their method shows that a modular form which vanishes 
on the image of $\ccJ_5$ in $\ccS_5/\Gamma_5(4,8)$
must vanish to order at least two on the image of $\ccJ_4$ 
in a boundary component of $\ccS_5/\Gamma_5(4,8)$. Since our modular forms on $\Gamma_5(4,8)$ have weight $8$ and also the defining modular form for $\ccJ_4$ has weight $8$, these modular forms must be identically zero on each boundary component of $\ccS_5/\Gamma_5(4,8)$. See also \cite{FO} for cusp forms of low weight on $\ccS_5$.

For recent progress on the characterization of $\ccJ_5$ in $\ccS_5$ using Schottky-Jung relations and geometrical methods, we refer to \cite{S}.

\

\end{document}